Andrei Rodin

Ecole Normale Supérieure

Towards a Hermeneutic Categorical Mathematics

or

why Category theory goes beyond Mathematical Structuralism

1.<u>Mathematical Interpretation</u>

*a) Hermeneutics of the Pythagorean theorem*

In his (1883) and a number of later writings W. Dilthey urges the autonomy of humanities (Geisteswissenschaften) from natural sciences (Naturswissenschaften). Dilthey's principle argument is methodological: objecting against attempts of Compte and Mill to extend scientific methods to moral, political and other humanities issues Dilthey purports to constitute the autonomy of humanities by providing them with a proper methodology independent from that of natural sciences. An important role in the Dilthey's methodology of humanities plays his notion of *hermeneutic understanding,* i.e. an understanding achieved through an interpretation, which Dilthey distinguishes from understanding achieved through a scientific explanation. As it has been noticed already by Husserl (1954) and recently stressed by Brown (1991), Crease (1997) and Salanskis (1991) hermeneutic issues are in fact not less important in natural sciences and mathematics than in humanities. A straightforward evidence of the relevance of understanding through interpretation in mathematics comes from the usual school practice: what counts as a genuine understanding (as opposed to mechanical memorising) of a given mathematical fact by a pupil is his or her capacity to formulate and prove it in his or her "own words" and apply it in a new unexpected situation. Obviously in a research environment the variability of forms of expressions of mathematical contents is even higher.

How this variability allows for a stable translatable mathematical content (meaning) and how precisely this phenomenon can be described? I think that this question has been so far very little studied. The question is not really specific for mathematics and can be easily reformulated as a problem of general theory of meaning. However the case of mathematical meaning is particularly important because in this case one can easier see how to treat the problem by rigorous mathematical methods. As far as the problem is solved inside mathematics one can think about application of the obtained mathematical solution elsewhere.

It might be argued that I confuse here a mathematical content with cognitive and social activities through which this content is "proceeded", or with symbolic and linguistic *forms* in which this content is expressed and communicated. It might be further argued that once a mathematical content is considered on its own rights hermeneutic issues become irrelevant. In the next paragraph I shall show that the latter claim is wrong because there are issues, which can be described both as "purely mathematical" and hermeneutic. In the present paragraph I shall show that the former claim is problematic (even if not plainly false) because a mathematical content cannot be easily separated from its symbolic form in a way making this form mathematically irrelevant.

Consider Pythagorean theorem. As formulated in (Lang&Murrow 1997, p.95) the theorem says this:

(LM) Let $XYZ$ be a right triangle with legs of lengths $x$ and $y$, and hypotenuse of length $z$. Then $x^2 + y^2 = z^2$.

(Doneddu 1965, p.209) under the title *Pythagorean theorem* states the following (my translation from French):

(D) Two non-zero vectors $x$ and $y$ are orthogonal if and only if $(y-x)^2 = y^2 + x^2$.

Finally the famous proposition 47 of Book 1 of Euclid's *Elements* states this (hereafter I quote Euclid by Heath's translation (1926):

(E) In right-angled triangles the square on the side subtending the right angle is equal to the squares on the sides containing the right angle.

Are (LM), (D) and (E) different forms of the same theorem? An attentive reader will immediately answer in negative pointing to the fact that (D) comprises the proposition usually called the *converse* of the Pythagorean theorem. This Donnedu's terminological decision I leave without commenting but ask the same question exchanging (D) for its *only if* part (which I denote (D') for further references). Now a plausible answer is *Yes, of course!* However obvious might be the answer, particularly in a mathematician's eyes, let's make some basic hermeneutics about (LM), (D'), (E) and read these propositions carefully before making a decision. Obviously each of the three propositions can be correctly interpreted only within a larger theory. Fortunately all the three books, from which I took the quotations, are elementary textbooks and so require no or very little previous mathematical knowledge. So in each of the three cases it is quite clear what is the corresponding larger theory. To simplify my task I shall leave out almost everything concerning *proofs* of the theorems and discuss only their statements.

Lang&Murrow in their book for beginners provide a fairly minimalist conceptual basis for their version of the Pythagorean theorem: before learning the theorem a student is supposed only to habituate him- or herself to basic geometrical constructions like triangles and learn the notion of *length* of a given straight segment. The latter notion reduces in the Lang&Murrow's book to the notion of *distance* between two given points. Lang&Murrow introduce the notion of distance through informally stated axioms of metric space and occasionally mention that distances are *numbers* one reads off from a graduated ruler. What a smart kid should think given two different rulers one of which is graduated in inches and the other in centimetres? The unwillingness of the authors to elaborate on this point is understandable since they commit themselves to keeping the Pythagorean Secret (not to be confused with the Pythagorean theorem) out of the reach of their students. I mean the incommensurability problem. Following the teaching strategy, which the legend attributes to Pythagoras himself, Lang&Murrow like most of authors of today's elementary mathematical textbooks reserve the truth about incommensurability for those of older students who choose to study mathematics at an advanced level and are expected not only to see the problem but also treat it by a modern remedy.

(Doneddu 1965) shows how the remedy may look like. This textbook applies Bourbaki's "architectonic" principles: it starts with making up a Boolean set-theoretic framework, then develops on this basis a theory of real numbers, and only after that comes to

geometrical issues. Euclidean space is construed here as a vector space over the field of the reals. Formula

$$(y - x)^2 = y^2 + x^2$$

requires an accurate interpretation: the minus sign on the left denotes the subtraction of *vectors* while the plus sign on the right denotes the sum of real numbers (so the two signs do *not* denote here reciprocal operations); the squares on both sides stand for the scalar product of vectors. We see that the price of the rigor is quite high: Donnedu's version of Pythagorean theorem requires a much more serious preparatory work, and after all the theorem doesn't explicitly refer to any triangle at all!

Euclid's classical presentation of the Pythagorean theorem depends, of course, on principles laid out in the beginning of the first book of the *Elements*. An extended historical comment on this theorem wouldn't be appropriate here, so I shall stress only one point often presenting a difficulty for a modern reader. (E) says that the bigger square equals to the smaller squares. How to understand this? An interpretation that immediately comes to mind is this: the *area* of the bigger square equals the *sum of the areas* of the smaller squares. But this is certainly *not* what Euclid says. Euclid speaks here about equality of figures, not about equality of their areas. In spite of differences between Euclid's and Hilbert's axiomatic methods (stressed in what follows), it is helpful to think about Euclid's notion of equality as formally introduced through Axioms of the first Book of the *Elements*. The Axioms give this: (i) equality is symmetric and transitive (Axiom 1; symmetry of this relation is granted by the linguistic form, in which the Axioms are expressed); (ii) congruent figures are equal (Axiom 4); (iii) figures composable of equal figures or complementable to equal figures, are equal (Axioms 2,3). This accounts for Euclid's equality as a binary relation, so (E) saying that one square equals to two other squares still remains puzzling. The obvious solution of the puzzle is this: the union of the two smaller squares is one relatum here while the bigger square is the other. Applying (iii) to the (topologically disconnected) union of two squares and using other Axioms and a number of preceding theorems Euclid proves his version of Pythagorean theorem. The equality of areas, of course, implies the Euclidean equality just explained. The converse, however, doesn't always hold. It doesn't hold, for example, for circles: given two smaller circles such that the sum of their areas is equal to the area of a bigger circle the bigger circle is obviously not equal in the Euclidean sense to the two smaller circles. A less trivial mathematical fact is that polyhedra of the equal volume are, generally speaking, not equal in the Euclid's sense either.

This short exposition of (LM), (D') and (E) is by far sufficient for claiming the obvious: although there is a sense in which all the three propositions express the "same theorem", the differences between them are anything but merely linguistic. So the claim that the three propositions say "essentially the same thing" shouldn't be taken at its face value. One who seeks to sweep the issue of interpretation out of mathematics might take a different strategy and argue that, say, only (D') represents the Pythagorean theorem in its correct form while (E) is hopelessly outdated and (LM) is just a simplified account for kids. Then it may be argued that the problem of translation between (LM), (D') and (E) belongs to history of mathematics and to mathematics education rather than to pure mathematics. It is obvious however that (LM), (D') and (E) are called by the same name of the Pythagorean theorem not for purely historical reason but because of similarity of their mathematical contents. How to grasp this similarity in rigorous terms? How the claim that given mathematical propositions A, B "say essentially the same thing" can be possibly justified?

In order to answer these questions consider first the notion of logical equivalence. May one generally look at logically equivalent mathematical propositions A, B (i.e. propositions which imply each other) as different expressions of the same mathematical fact? Obviously not. Equivalent mathematical propositions may "mean" very different things, so the equivalence of two given mathematical propositions may be very far from being obvious. For example in the traditional Euclidean setting the theorem saying that the sum of internal angles of any triangle is equal to two right angles is equivalent to the to the Fifth Postulate ("Axiom of Parallels").  But this cannot be seen immediately without a proof because the two propositions *mean* quite different things. So the logical equivalence is not sufficient for the purpose. But is it necessary? Given that A, B are different formulations of the "same theorem" is it always the case that A, B are logically equivalent (in symbols A⇔B)? This might sound like a reasonable requirement but looking at our example of Pythagorean theorem one would wish to relax it. For one cannot assert and moreover prove A⇔B unless A, B belong to the same theory. But our

(LM), (D'), (E) all belong to different theories! Let's see whether the wanted background theory C can be acquired. Obviously C should incorporate not only propositions (LM), (D'), (E) but also their corresponding mother theories or at least relevant fragments of these theories. Although some *bricolage* of this sort can be certainly made (school geometry textbooks provide many such examples) a coherent theory systematically combining Bourbaki's setting with Euclid's (in its original form) can be hardly imagined. Instead of combining theories treating differently the "same subject" (in a sense that we are still looking to define) one may rather *interpret* these theories in each other's term. Let's see how this works in the given example.

(E) translates into (LM) through replacement of squares constructed on sides of a given triangle by squares of lengths of these sides, and replacement of Euclidean equality explained above by equality of real numbers. As far as (E) and (LM) are considered in isolation from their mother theories this translation may be viewed as reversible. Clearly this is the existence of such translations (I leave translations of (D') into (E) and (LM) to the reader) that allows one to speak of similarity or identity of contents of different formulations of Pythagorean theorem. As far as these translations are well accounted for the problem of identity of contents or "meanings" can be sorted out through a reasonable convention. However such a convention cannot be reasonably made ad hoc, so one needs a general theory of translation between mathematical theories.

One might argue that formal axiomatic method reinforced by symbolic logic and Model theory provide what is wanted here. This method relies onto an assumption that theories allowing for mutual term-wise reversible translations are essentially the same and can be reduced to a single "formal" theory. So the formal method provides at least a sufficient condition of "essential sameness" of mathematical theories. However this condition is too strong to be useful in the case of Pythagorean theorem and other similar cases. Although the aforementioned method of translation of (E) into (LM) can be applied in wider contexts (in particular to propositions of Book 2 of *Elements*) it doesn't apply to the whole of *Elements* and doesn't respect its deductive structure. So the suggested translation doesn't make Euclid's and Lang&Murrow's accounts into two isomorphic models of the same formal theory. As far as mother theories of (E) and (LM) are taken into consideration (as they should be) the suggested translation doesn't look like isomorphism any longer. Although Model theory allows for more involved relationships between models than isomorphism the very notion of model assumes a notion of formal theory (that a model in question is a model of). But as our example clearly demonstrates one can have a sound notion of translation between non-formalised theories. Moreover it seems reasonable to consider formalisation as a translation of a special kind. So formal method and Model theory can hardly replace a theory of mathematical translation (which can be also called a theory of mathematical hermeneutic) we are looking for.
In what follows I shall argue that the wanted theory is Category theory.

*b) The rise of interpretative mathematics*

Hermeneutic issues about mathematics discussed in the previous paragraph involved historical and educational dimensions. Now I shall show how the notion of interpretation got involved in mathematics in a more abstract manner and so became "purely mathematical". The term "interpretation" appears in the title of Beltrami's paper of 1868 *Saggio di interpetrazione della geometria non-euclidea* which in eyes of many first showed that non-Euclidean geometry was something real. However it seems appropriate to start the history of interpretation as a mathematical concept not with this Beltrami's paper but from Gauss' geometrical work of 1820-ies (making a reasonable distinction between history and corresponding prehistory). Interestingly (but perhaps not so surprisingly) interpretation became a genuine mathematical issue during the same period of time and in the same part of Europe where and when Schleiermacher, Dilthey and their followers stressed the role of interpretation in humanities. Let me recall the story. Lobachevsky discovered the non-Euclidean geometry presently called by his name through "playing with axioms", namely through replacement of Euclid's Fifth Postulate ("Axiom of Parallels") by its negation. Like his numerous Ancient and Modern predecessors Lobachevsky hoped to get a contradiction and hence a proof of the Postulate. However like Bolyai and few other people working on the problem around the same time Lobachevsky at certain point changed his attitude and came to the conviction that he explored a new vast territory rather than approached the desired dead end. He called this new geometry *Imaginary* (Lobachevsky 1837) because of the speculative character of his enterprise and probably as a precaution: if his theory would turn after all to be contradictory he would win anyway getting the wanted proof.

While Lobachevsky made his discoveries following a traditional line of research Gauss got a totally different insight on the problem. Although Gauss' name hardly needs an additional promotion in the history of mathematics, I claim that his role in the discovery of non-Euclidean geometries is often misinterpreted and underestimated (like in Bonola 1908). After reading Bolyai's paper (1832) Gauss claimed that he found for himself nothing new in it, and this claim is at least partly confirmed by existing evidences. He made then a controversial remark that "to prise it (the Bolyai's paper - *AR*) would mean to prise himself" (Bonola 1908). Historians often explain this Gauss' reluctance by personal and sociological reasons or by his alleged epistemological conservatism. I think that Gauss' cautious attitude to Bolyai's and Lobachevsky's results doesn't need such non-mathematical explanations. Gauss didn't accept the notion of geometry as a speculation but considered it as an empirical science. For this reason he was very sceptical about the whole line of research that led to Lobachevsky's discoveries (known at the time as the "theory of parallels"). This Gauss' attitude had nothing to do with conservatism: on the contrary, in Gauss' eyes the theory of parallels was too *traditional* and missed really new ideas[1]. From the today's historical distance it is easier to argue that Gauss was perfectly right, and that his insights later spelled out by Riemann played more important role in the change of views on space it time occurred in the 19-20-th century than the whole story about the Fifth Postulate. However one had to be a mathematician of the rang of H.Weyl to see this clearly already in 1918:

"The question of the validity of the "fifth postulate", on which historical development started its attack on Euclid, seems to us nowadays to be a somewhat accidental point of departure. The knowledge that was necessary to take us beyond the Euclidean view was, in our opinion, revealed by Riemann."

Given that Riemann's concept of manifold provides the mathematical basis of the today's best theory (or theories) of space and time, and so replaces in this role Euclidean space of Classical mechanics, Weyl's point is hardly disputable. That Riemann's geometrical works are directly based on Gauss' is not disputable either. What makes it difficult for a part of historians to appreciate Gauss' and Riemann's contribution is apparently the fact that the mathematical work of these people doesn't fit the popular story about liberalisation of mathematical thought from its alleged stickiness to everyday spatial experience by Lobachevsky. Riemannian geometry just like Euclidean geometry about two and a half millennia earlier has been first sketched on the ground by Gauss and only after that worked out in a more abstract form by Riemann and successfully applied by Einstein to Heavens. What triggered this whole development was a new attentive look at the space we live in rather than a mere play of imagination or an abstract mathematical speculation.

---

[1] « In der Theorie der Parallellinien sind wir jetzt noch nicht weiter als Euklid war. Diess ist die partie honteuse der Mathematik die frue oder spaet eine ganz andere Gestalt becommen muss .» *Werke* v.8, 166. This is written in 1813, that is, before Bolyai's and Lobachevsky's works were published. But I believe that Gauss didn't find the new *Gestalt* he was looking for neither in Bilyai nor in Lobachevsky.

In 1818-1832 Gauss was busy with what geometry used to be in its early age and later got a different name of *geodesy*. He started with the obvious observation that the hilly terrain of Hanover was not Euclidean plane. He also saw that that the current physical hypothesis according to which the Kingdom of Hanover together with its mother planet float in the infinite Euclidean space was not particularly helpful for geodesic purposes. So he looked for different geometrical models. This led him to the theory of curved surfaces that he presented in his *Disquitiones generales circa superficies curva* published in 1827. There are firm evidences that Gauss saw connections between non-Euclidean (or *anti-Euclidean* as he himself called it) geometry obtained through playing with axioms and the geometry of curved surfaces he was working on. Gauss rightly guessed that the latter leads to a more fundamental generalisation of the notion of space than the former.

The key Gauss' idea that allowed for this generalisation was the idea of *intrinsic geometry* of a given surface. Abbott's popular *Flatland* (first edition 1884) explains the idea but oversimplifies the general situation: really interesting things happen when one considers living on a curved surface rather than in the Flatland, for example on a sphere like our globe. Abbott seems to suggest to the reader the following moral: just like 3D creatures like ourselves are in a position to observe things going on a plane from a "higher viewpoint" and perform tasks impossible on the plane (like escaping from a plane prison) a creature living in a space of 4 or more dimensions would find herself in a similar position with respect to us ordinary humans. Abbott might believe that this higher viewpoint could be achieved through doing mathematics. But this moral is hardly justified mathematically. However fascinating the idea of 4D space might be the intuition that rising of dimension always allows for solving problems in lower dimensions is quite misleading. Given a geometrical problem on a plane switching to 3D space is rarely helpful. The idea of "intrinsic viewpoint" obtained through *lowering* the dimension is much more profound, and as a matter of fact it played a far more important role in mathematics of late 19-th and the whole of 20-th century.

The possibility of purely intrinsic description of a surface leads to a generalisation of the Euclidean notion of space, that is, to non-Euclidean geometries: only in the special case when the given surface is flat its intrinsic geometry is Euclidean. Riemann in (1854) sketched this new notion of space and called the new concept by the term *manifold* occasionally used before by Gauss. A Riemannian manifold is a $n$-dimensional analogue of a curved surface seen intrinsically. The talk of a *curved space*, which became colloquial after Einstein, refers to the concept of Riemannian manifold.

Let's now return to Beltrami. This man like Gauss started his research in geometry with geodesy. He knew Gauss' results in this domain and tried to elaborate on them. Beltrami's geometrical discoveries originated from the classical cartographic problem: How to make a plane map of a curved surface? More specifically Beltrami asked the following question: How to map a curved surface onto a plane in such a way that the mapping is one-to-one on points and *geodesic* lines on the surface go to straight lines on the plane? (A *geodesic* is a line that marks the shortest path between its close points; geodesics on a plane are straight lines. So the notion of geodesic generalises upon that of a straight line for the case of curved surfaces. The notion of geodesic is *intrinsic*: distances of paths between points of a surface don't depend on how the surface is embedded into an outer space.) Here is first important result obtained by Beltrami: such a mapping is not possible unless the *curvature* of the surface in question is constant. (The *curvature* is a basic *local* intrinsic property of a given surface, which shows how much the surface is curved around a given point; the concept is due to Gauss.) Hence Beltrami's interest to surfaces of constant curvature (a sphere is an obvious non-trivial example) - the issue which had been already studied before Beltrami by another Gauss' follower Minding. In 1866 Beltrami read (Lobachevsky 1840) in French translation recently published by Beltrami's long-term collaborator Houel. This led Beltrami to his main discovery presented in his *Saggio*: mapping geodesics of a surface of constant *negative* curvature (which Beltrami called a *pseudo-sphere*) to straight lines of a plane one gets the Lobachevskian "imaginary" but not the "real" Euclidean geometry. Beltrami's understanding of this result was this: Lobachevsky's "imaginary plane" is *in fact* nothing but a pseudo-sphere!

*Saggio* was impressive but had two serious flaws. First, Beltrami didn't notice that his model was only partial: it allowed to think of finite segments of geodesics of pseudo-sphere as of segments Lobachevskian straight lines but didn't represent the whole lines. This fact has been noticed by Helmholz and Klein soon after Beltrami's publication. In (1901) Hilbert showed that the difficulty couldn't be overcame through replacing the pseudo-sphere by another surface. The other flaw of the *Saggio* is stressed by Beltrami's himself in the end of this work: the suggested interpretation of Lobachevskian planimetry doesn't generalise to the 3D case. For that reason Beltrami at certain point called the Lobachevskian stereometry a "geometrical hallucination". The need of very different treatment of 2D and 3D cases made Beltrami to suspect that something went wrong. After reading Riemann's *Habilitaetsvortrag* Beltrami found what he now thought was a better answer to the question "What is Lobachevskian space?": this is *in fact* a Riemannian manifold of constant negative curvature. This latter answer holds for any number of dimensions.

From an epistemological viewpoint the story just told is very unusual. A classical model of science outlined in Aristotle's *Analytica Posteriora* requires first to fix a subject matter of study, then study it and organise results of the study into a theory. It is assumed that the study reveals earlier unknown features of the subject matter but doesn't replace it by something else. In the case of mathematics this subject matter could be understood as *figure*, *number*, *conic sections*, *space* or something more specific. But the Lobachevsky-Beltrami's case turned this traditional epistemological scheme upside down. Lobachevsky came up with an "imaginary" mathematical theory, which didn't have any clear subject matter but was too good to be simply abandoned. Beltrami looked for an appropriate subject matter for this theory and found at second attempt a satisfactory answer. At this point the talk of interpretation apparently became redundant and the story about a theory looking for its subject could be classified as a historical curiosity. However at the time there were other reasons to take the idea of freedom of interpretation in mathematics seriously. An earlier example of interpretation, which Beltrami mentions in Saggio, is the representation of complex numbers as points of Euclidean plane. Another important example, which later motivated Hilbert, is the duality principle in Projective geometry. The history of this latter geometrical discipline is closely connected to the history of non-Euclidean geometries outlined above but I shall not explore it here and only briefly mention the history of the duality principle. As it has been first noticed by Poncelet in 1822 given a theorem of Projective geometry one may get another theorem formally exchanging words "points" and words "straight lines". So points can be "interpreted as" straight lines and straight lines as points without changing the theory. Gergonne called this phenomenon *duality* and Steiner made it into foundations of Projective geometry (for a detailed historical account see Kline 1972, pp. 845-46).

In spite of its striking new features revealed in 19-th century the notion of mathematical interpretation was not entirely new. One may even argue that it was known in mathematics since its early history. For centuries mathematicians used to represent some objects by other objects and substitute symbols for other symbols looking for invariants of such operations. This is what the whole discipline of *algebra* is about: algebraic *variables* take different values leaving the *form* of a given algebraic expression invariant. In Antiquity mathematicians used to represent (natural) numbers by geometrical objects (think about Pythagorean *figured numbers*, see Heath 1965, v.1) and in Modern times they learnt to represent geometrical objects (points) by numbers. In both cases it is assumed that a given represented object remains invariant under exchanges of representations. It may be further argued that the idea of invariance through representations is fundamental for the very notion of mathematical object: for example, a *circle* (a mathematical object) may be thought of as an invariant of any series of exchanges of material objects (drawings and the like), which are told to *represent* the mathematical circle. This view has been thoroughly spelled out by Plato who in fact suggested a more elaborated theory according to which mathematical objects in their turn represent (are *images of*) things of yet another sort he called *ideas*. One may also tell a plausible story about representation of counted material objects (say, caws) first by other material objects like pebbles, which are better manageable, then by written symbols and finally by abstract mathematical units. Although this story assumes the order of representation which is opposite to Plato's, the idea that the interface between the material and the mathematical worlds can be described in terms of representation remains unchanged.

Let me now stress the principle point of this paper. Although the older notions of representation and substitution (and the related notion of invariance through substitution), on the one hand, and the notion of *interpretation* as it emerged in the geometry of the late 19-th century, on the other hand, had indeed much in common the latter does *not* reduce to the former (while the former indeed reduces to the later). I claim that this fact has been only partly recognised in the end of 19-th - the beginning of 20-th century by Hilbert and other people who treated the new situation along the traditional pattern mentioned above. So my task is now to explore possibilities left out by this mainstream development.

2. Foundations

a) Elements *and* Grundlagen

We have seen that the non-Euclidean geometry emerged in 19-th century had two well distinguishable sources. The first is the traditional line of research aiming at proving the Fifth Postulate through drawing a contradiction from its negation. This line of research started in Antiquity and resulted into Bolyai and Lobachevsky's works. The second line started with Gauss' geodesic work and led to the notions of intrinsic geometry and Riemannian manifold. Beltrami brought the two lines together showing that Lobachevskian spaces *are* Riemannian manifolds of a particular sort. However the question of *foundations* of the new geometry remained open. Obviously neither Euclid's *Elements* in its original version nor the generalised version of Euclid's system proposed by Bolyai (his *absolute geometry*) could serve this purpose. Taken seriously the problem of foundations of geometry in the end of 19-th century would have to account not only for the Riemannian geometry in its full generality but also (at least) for projective geometry and topology. Klein (1893) made a substantial progress toward a theoretical unification of geometry developing various links between these disciplines but he didn't produce anything like a replacement of Euclid's *Elements*. Hilbert's *Grundlagen* first published in 1899 partly meets this challenge.

I say "partly" because in this work Hilbert accounts only for a very limited part of his contemporary geometry. Basically the *Grundlagen* shows how the Euclidean geometry looks like in a new context including Lobachevskian geometry and some other geometries obtained through "playing with axioms" of Euclidean geometry. *Grundlagen* provides an effective framework and non-trivial examples of such a logico-mathematical game and treats related logical issues of consistency and independence of axioms. So the *Grundlagen* continues the traditional Euclidean-Lobachevskian line and doesn't touch upon the Riemannian viewpoint. Hilbert's work became highly influential because of its *method*, viz. Hilbert's *axiomatic method*, not because of its content. Hilbert believed that using this method one might build appropriate foundations of the whole of mathematics and of other sciences. Veblen (1904) and many other people enthusiastically followed this idea. (For the history of Hilbert's *Grundlagen* see Toepell 1986).

As it is often happens with projects aiming at reform of the whole system of human knowledge Hilbert's project of axiomatisation of mathematics and sciences brought controversial results. On the one hand, nothing like an effective global axiomatisation of mathematics, and moreover of natural sciences, has been ever achieved. On the other hand, Hilbert's *Grundlagen* remains a paradigm of a "reasonably formal" (as opposed to "purely formal") axiomatic system in eyes of the majority of working mathematicians. A today's student of mathematics may easily think - and read in many textbooks - that the axiomatic method as it is presented in the *Grundlagen* is just a more rigorous version of the method first used by Euclid in his *Elements*. In this paragraph I shall try to show that this view is misleading, and stress a specific character of Hilbert's axiomatic method against that of Euclid. For this end let me briefly compare few first pages of the *Elements* and the *Grundlagen*.

Euclid starts with his famous definition of point while Hilbert assumes primitive notions of point, straight line and plane without trying to define them. There are two closely related but still different reasons why Hilbert does this. The first is Hilbert's epistemological holism: he assumes that (1a) unless a theory (or at least its axiomatic basis) is wholly given these basic concepts cannot be adequately conceived and moreover that (1b) an adequate conception of these concepts can be provided by axioms alone (without help of definitions). My point is that this holistic view is not incompatible with Euclid's *Elements*, and so not specific for Hilbert's approach. There is no reason to assume (as apparently did Hilbert) that Euclid's definition of point as "having no parts" is supposed to provide a full grasp of the concept. Euclid might grant (1a) and think of his definition of point as nothing but useful hint. (This interpretation is supported by what Aristotle says about definitions in his *Analytica Posteriora*.) Although (1b) is far more controversial Euclid might still agree that there is a special sense in which the definition of point given in *Elements* is redundant.

The other reason why Hilbert leaves points and lines undefined is more involved. It becomes clear from the following passage where Hilbert explains his method in the nutshell to Frege:

"... surely it is self-evident that every theory is merely a framework or schema of concepts together with their necessary relations to one another, and that basic elements can be construed as one pleases. If I think of my points as some system or other of things, e.g. the system of love, law and chimney sweeps ... and then conceive of all my axioms as relations between these things, then my theorems, e.g. the Pythagorean one, will hold of these things as well. In other words, each and every theory can always be applied to infinitely many systems of basic elements. For one merely has to apply a univocal and *reversible* one-to-one transformation and stipulate that the axioms for the transformed things be correspondingly similar. Indeed this is frequently applied, for example in the principle of duality, etc." (cit. by Frege 1971, p.13, italic mine).

Since a point is allowed to "be" (or "thought of") a "system of love or of chimney sweeps" (or a beer mug according to another popular Hilbert's saying) nothing like Euclid's definition of point can be any longer useful. Neither Euclid nor any other mathematician before Hilbert could agree that in geometry "basic elements can be construed as one pleases"! (Frege didn't agree with this either.) Epistemological holism (however radical) doesn't imply this view. But let's see a true mathematical reason behind Hilbert's colourful rhetoric. Notice Hilbert's reference to the duality principle in

Projective geometry. What Hilbert aims at here is the following: to construe a mathematical theory leaving its interpretation free, that is, to construe it "up to interpretation".

This also explains the second difference, which hits the eye when one compares first pages of *Elements* and *Grundlagen*. After giving basic definitions and before coming to axioms Euclid lists five *Postulates* while in Hilbert's work there are no such things at all. Although the *Grundlagen* is not the first introductory text in geometry written after *Elements* without making use of postulates the absence of postulates or their analogues in *Grundlagen* is remarkable. Let me now explain what are Euclidean postulates. Consider the first three Postulates of the *Elements*:

> **1.** *to draw a straight line from any point to any point*
> **2.** *produce a finite straight line continuously in a straight line*
> **3.** *describe a circle with any centre and distance*

and observe that unlike axioms Postulates 1-3 are *not* propositions about geometrical objects but descriptions of certain *operations* performed with geometrical objects. Unlike *Grundlagen Elements* is a system of *constructions* generated by a set of elementary operations described by the Postulates 1-3 (that is, constructions by the ruler and the compasses) *and* a system of propositions associated with these constructions. (Proclus in his *Commentarium* analyses the distinction between the two aspects of the theory of the *Elements* in terms of the Platonic ontological distinction between Becoming and Being: geometrical objects are treated by Euclid both *qua* constructed (generated) and *qua* pre-existing entities. But one doesn't need to buy the Platonic metaphysics to recognise the distinction.) So the view on *Elements* as a system of propositions obtained through a logical inference from a system of basic propositions called axioms is plainly wrong. Whether or not the theory of *Elements* can be reasonably reformulated as such a system of propositions is a different issue, which I leave aside. *Elements* doesn't fit the above description not because it is somewhat "imperfect" but simply because Euclid's project is different. But *Grundlagen* is supposed to fit this description indeed. What in the case of *Grundlagen* replaces postulates in their function of giving a "flesh" to geometrical construction is the assumption that, anachronistically speaking, the theory of *Grundlagen* may have a *model*. However *Grundlagen* doesn't provide by itself any special means of building such a model. Nothing like Euclid's Postulates 1-3 could be appropriate in *Grundlagen* where the interpretation of a given geometrical theory is supposed to be left free.

We see that Hilbert's idea of construing mathematical theories "up to interpretation" meets indeed in a certain way the challenge of "interpretative mathematics" of 19-th century. Let me now explain why I consider Hilbert's approach as traditional in spite of its strikingly new features. The subject matter of the traditional mathematics can be described in an Aristotelian vein through distinguishing certain properties of material objects called *mathematical* properties (like *shape*) and ruling out non-mathematical properties like *colour*. A mathematician is allowed to use material objects in his or her work, and even make some mathematical use of non-mathematical properties (think about the problem of four colours) but he or she should never confuse these material objects with the proper subject matter of a mathematical study. The new kind of mathematics invented by Hilbert (called *formal* mathematics) makes a further step in the same direction: it rules out *all* non-relational properties as irrelevant and allows into its proper subject-matter only *bare things* and *bare relations* between these things. In practice a mathematician may think about these new mathematical things in the usual way, call them by usual names and use convenient helpful drawings. But one is also left free to use some unusual names and images for it. This becomes is a matter of personal taste or, perhaps, of a research skill. In any event names and images don't count in the final result: a ready-made mathematical theory must not depend on traditional mathematical notions and on usual intuitions associated with these notions (let alone names and pictures) just like it must not depend on the colour of inks used for writing mathematical papers. Hilbert's "formal" approach assumes that the *same* mathematical theory can be interpreted through traditional constructions in different ways just like it can be written down by different inks. The *Grundlagen* provide a tentative theory (in fact few different theories) with desired unusual properties. Remarkably the subsequent development of Hilbert's idea of formal theory at a certain point allowed for a new definite answer to question What is subject matter of mathematics? The answer is this:

mathematics is about *sets*. It might seem that the key idea of Hilbert's project of "leaving interpretation free" rules out a possibility of such a definite answer however this is not the case. The new view on the subject matter of mathematics connects to Hilbert's project roughly as follows. Hilbert proposes to think about geometrical points and straight lines (or other primitive mathematical notions) as abstract "things" forgetting about their convenient images. That things in question form *sets* appears to be a natural and technically speaking very useful assumption. Hence the view on Set theory as a "general theory of things" relevant to mathematical matters. (The distinction between "things" or "urelements" and sets of such things assumed by Cantor has been later largely abandoned by Zermelo's suggestion as redundant from a formal viewpoint.) From a historical perspective Hilbert's project can be then viewed as aiming at a higher level of mathematical abstraction and bringing about a new kind of abstract objects, namely abstract sets. It goes without saying that Hilbert is not the only person to be credited for this achievement. Anyway the scandalous situation of the end of 19th century when mathematics seemed to have loosen its subject matter has been successfully resolved: with the notion of abstract set mathematics regained a firm ground.

Let me now show that this solution is not satisfactory and suggest a remedy. I shall not speak about internal problems of set-theoretic foundations widely discussed in the earlier literature but challenge the general view according to which a formal foundations is appropriate for mathematics.

b) *Categoricity and reversibility*

Suppose a formal system S of (uninterpreted) axioms has two different models A, B. The formalist viewpoint outlined above suggests to disregard differences between A and B as mathematically irrelevant like the difference of colour of two drawn circles. But suppose that now the system S is extended by some additional axioms, and that A is a model of the extended system S' but B is not. Since the difference between A and B is now grasped by the formal method a formalist must recognise the difference between A and B as essential. So in order to be consistent a formalist needs a criterion of the "essential sameness" of models independent of their corresponding theories. In the letter quoted above Hilbert says clearly what criterion of "essential sameness" of models he has in mind: this is isomorphism, i.e. existence of reversible element-wise translation between the models. In other words his idea of formal theory is that of a theory construed up to isomorphism (of interpretations). But what guaranties that a given theory like that of *Grundlagen* is indeed formal in this sense, i.e. that all its models are in fact isomorphic? The desired property has been called by Veblen (1904) *categoricity*. Beware that this term has nothing to do with Category theory.

When Hilbert first published his *Grundlagen* he apparently didn't ask this question. He did this in the period between publications of the first and the second editions of *Grundlagen* (1899-1903). So in the second edition of *Grundlagen* published in 1903 Hilbert tried to safeguard categoricity of his theory by a dubious axiom postulating that any model of his theory is maximal in the sense that it is impossible to obtain another model of the theory through introducing new elements into it. The existence of such a maximum seems to be incompatible with the standard model theory.

In 20-th century people learnt to be more tolerant to the lack of categoricity. Remind that Zermelo-Frenskel axiomatic Set theory (ZF), Peano Arithmetic (PA) and some other theories commonly viewed as important turned to be non-categorical. To preclude the right of these theories to be qualified as formal on this ground would apparently mean to go too far. To save the situation philosophers invented the notion of "intended model", that is of model chosen among others on an intuitive basis. Isn't this ironic that such a blunt appeal to intuition is made in the core of formal axiomatic method? I agree with F. Davey who recently argued that "no-one has ever been able to explain exactly what they mean by intended model". (FOM, 13 Oct 2006). Other people question the categoricity requirement. R. Lindauer: "Why rule out non-standard models of 1st-order PA? What's wrong with having other models? Why should we be making our model-world smaller and not larger?" (FOM, 27 Oct 2006). Myself I believe that the lack of categoricity of theories like ZF and PA is indeed a serious flaw. At the same time I agree with Lindauer and other people who think that the pursuit of categoricity is misleading. These two claims might seem to contradict each other but they don't. Instead of pursuing categoricity or looking for a philosophical excuse of the lack of categoricity of popular formal theories I shall suggest in what follows to use an alternative method of theory-building.

I don't consider pointing to the problem of categoricity as a decisive argument against the formal method. I think that this problem is in fact quite artificial and originates from what I consider to be the principle flaw of the formal method. The flaw is the following. Remind the problem Hilbert faced suggesting his formal method: How to formulate a mathematical theory leaving its interpretation free? or How to formulate a theory "up to interpretation"? Here is Hilbert's response: such a theory **T** must be *formal*, which means that its primitive terms (objects and relations) are *variables* taking their semantic values ("meanings") through interpretations; given such interpretation (model) **M** one obtains another model **M'** of the same theory through a one-to-one substitution (exchange) of primitive terms. If **T** is categorical then the substitution of terms allows one to obtain *all* models of this theory from any given model. We see that Hilbert considers *reversible* transformations between models (one-to-one substitutions of terms) as the *only* kind of interpretation he has to cope with. But the notion of interpretation as it has emerged in geometry of the 19-th century does *not* reduce to such reversible interpretations (isomorphisms). Interpretations are, generally speaking, *non-reversible*. A theory determined "up to isomorphism" is not a theory determined "up to arbitrary interpretation". So Hilbert's formal method doesn't resolve the problem posed by the new "interpretative mathematics" but in a special case. The lack of a categoricity of workable formal systems is, in my view, nothing but a symptom of the latter problem.

To show that mutual interpretations of mathematical theories are, generally, non-reversible, let me come back to Beltrami. In terms of *Teoria* the principle result of *Saggio* can be formulated as follows: a 2-dimensional Riemannian manifold of constant negative curvature (= Lobachevskian plane) is embeddable (in fact only locally embeddable) into 3-dimensional Euclidean space which is another Riemannian manifold. The embedding is *not* reversible: one cannot embed the Euclidean space into the Lobachevskian plane. One may remark that we are talking here about a map between spaces (manifolds) but not about an interpretation between *theories*. But it is obvious that however the notion of theory is construed in this case the situation remains asymmetric: while Lobachevskian plane geometry can (modulo needed reservations) be explained in or "translated into" terms of Euclidean 3D geometry the converse is not the case.

An embedding is a sort of transformation called *monomorphism*. Carving a pseudosphere out of its ambient space one may show indeed that the pseudosphere is isomorphic (in fact only locally isomorphic) to the Lobachevskian plane. However this reasoning is misleading: in the given context a pseudosphere cannot be carved out from the Euclidean space and considered as a self-standing object. For if the pseudosphere is indeed carved out from its ambient space and considered as a self-standing space (manifold) it ceases to be Euclidean.

We see that the notion of interpretation in mathematics as it has emerged in geometry of the 19-th century doesn't reduce to the old idea of reversible term-wise *substitution* of values of variables, which has led Hilbert to his formal axiomatic method.

A similar point can be made about arithmetical models of plane Euclidean and other geometries used by Hilbert in *Grundlagen*. Perhaps one can indeed imagine geometrical points as usual dots, "systems of loves" or beer mugs indiscriminately. But representation of points by pairs of real numbers (or pairs of elements of another appropriate algebraic field) is a different matter. Unlike dots and beer mugs numbers are mathematical objects on their own rights belonging to a different mathematical theory, namely Arithmetic. "Translations" of geometrical theories into Arithmetic used by Hilbert are evidently non-reversible: they allow to recast geometrical theories in arithmetic terms but not arithmetic theories in geometrical terms. Hilbert certainly saw this. He didn't mean to say that Geometry and Arithmetic seen from a higher viewpoint turn to be the *same* theory (in fact he considered a possibility of reduction of one theory to the other but this is a different issue). Nevertheless he thought about this translation as an isomorphism, namely an isomorphism between basic geometrical objects and relations, one the one hand, and specially prepared arithmetical constructions, on the other hand. But such constructions obviously cannot be made outside an appropriate arithmetical theory! As far as this "target" arithmetical theory is wholly taken into consideration the translation in question doesn't look like an isomorphism any longer. We see that the example of projective duality mentioned by Hilbert (to leave alone the "system of love, law and chimney sweeps") is special and cannot be used for treating the notion of interpretation in mathematics in the general case.

Let me now suggest how the "hermeneutic challenge" can be more effectively met.

c) *Forms and Categories*

The concept of Riemannian manifold like many other mathematical concepts allows for a corresponding notion of map or transformation (between different items falling under the given concept). It is obtained as a straightforward generalisation of the notion of isomorphism of manifolds through giving up the reversibility condition. Maps between Riemannian manifolds are *differentiable* transformations. Think about these maps as usual maps in cartography. So they can be conceived as "interpretations" of manifolds in terms of other manifolds. Describing the concept of Riemannian manifold "up to interpretation" amounts to construing the manifolds together with all maps between them. A structure of this sort is called a *category*. The notion of category introduced by Eilenberg and MacLane in (1945) is very weak: one requires only an operation of *composition* of transformations between the objects, the associativity of this composition and the existence of an identity transformation for each object. Transformations of objects are called in the Category theory *morphisms*. These requirements are obviously satisfied by Riemannian manifolds and differentiable transformations. Other standard examples include sets and functions, groups and group homomorphism, topological spaces and continuous transformations, etc. A shall denote the category of Riemannian manifolds **RM** for further references.

The notion of category can be apparently thought of as a very general *form* of mathematical concepts. This might seem natural but in fact is quite misleading. For the notion of *form* can be described as invariant of a group of isomorphisms between objects of a given class which are told to have the "same form". More precisely a form is what one gets through identification of isomorphic objects. Think about Euclidean circles and usual geometrical transformations between the circles like motions and scalings. Such transformation are reversible, otherwise they wouldn't form a group. Changing the chosen class of transformations one changes the corresponding notion o form. For example through allowing for all reversible continuous transformations between circles one obtains the topological notion of circle (a "topological form"), which differs from the convenient metrical notion of circle (in particular the former unlike the later doesn't distinguish between circles and ovals). The suggested definition equally fits the colloquial notion of algebraic form: the isomorphisms associated with algebraic forms are substitutions of variables. It also fits the traditional Platonic conception of natural numbers as forms: think about bijections between finite sets.

Now if one think of the abstract notion of category, on the one hand, and independently construed mathematical objects (manifolds, groups, sets, etc) with their associated maps, on the other hand, one may think of substitutions of abstract categorial objects and morphisms by "concrete" manifolds and differentiable transformations, groups and group homomorphisms, sets and functions, etc. However these substitutions are not reversible because the categories in question are not isomorphic! This can be seen clearly through the standard construal of manifolds, groups, topological spaces, etc. as "structured sets"; in this case morphisms between manifolds, groups, etc. can be described as maps between base sets preserving corresponding structures. As far as these structures are different the morphisms are different too. So the abstract notion of category cannot be thought of in terms of reversible substitutions of abstract objects and abstract morphisms by some "concrete" terms. Hence an abstract category is not a form in the sense of the above definition.

Although different categories don't share anything like a common form they can be linked (or "transformed into each other" or "interpreted in each other's terms") by appropriate morphisms, which in this context are called *functors*. So different categories form various "categories of categories" just like different sets form sets of sets. One may even think (as does Lawvere in his 1966) about "the" category of categories albeit this risky notion is not necessary for my argument.

Observe that certain mathematical concepts like the concepts of circle, of singleton, of Euclidean plane and many others have this strong property: their instances are all isomorphic. I shall call such concepts *form-concepts*. But there is also a larger class of mathematical concepts having a weaker property: their instances make categories. The concepts of set, group, manifold, topological space and many others are not form-concepts but belong to this latter class. I shall call these latter concepts *category-concepts*. Form-concepts are category-concepts since a group of isomorphisms associated with a form-concept is a category having only one object and such that all its morphisms are reversible. But the converse is not the case. Thus a category is not a kind of form while a form (in the sense of the above definition) is indeed a kind of category. Remark that the concept of category is itself a category-concept while the concept of form is not a form concept. I shall not further elaborate here on the general concept of form but suggest that it is a category-concept as well.

The terminology just introduced helps me to formulate more clearly my objection to Hilbert's formal approach and make an alternative proposal. Hilbert's idea of "formal" theory amounts to conceiving of mathematical theories as form-concepts. (I'm talking about particular theories, not about the general concept of theory.) But the analysis given in the previous section suggests thinking about mathematical theories as category-concepts (since mutual interpretations between theories are generally non-reversible!). So my proposal amounts to conceiving of a mathematical theory as a category. No analogue of the categoricity problem (in Veblen's sense) appears within the category-theoretic approach, which I'm going now to explain.

d) *Categorical theories and functorial semantics*

The categorical method of theory-building can be roughly described as follows: theories and concepts are specified by their categorical properties (rather than by formal properties). The idea is this: given a category-concept (like that of manifold, group, etc.) to specify an "abstract" category which could be reasonably identified with a corresponding "concrete" category obtained through a standard (usually set-theoretic) construal of the concept in question and the corresponding notion of map. Observe that the usual distinction between abstract and concrete categories just mentioned should be abandoned as far as one takes a category-theoretic viewpoint to start with. As far as concepts of object and morphism are taken as basic and other mathematical concepts are supposed to be reconstructed in their terms there is no reason to call them "abstract".

For the case of sets such work has been first done in (Lawvere 1964) where the author has provided a list of axioms distinguishing the category of sets among other categories. The case of category of Riemannian manifold is much more involved, I refer the reader to (Kock 1981), (McLarty 1992) and (Bell, forthcoming) for details. Observe that categorical approach to building geometrical spaces generalises upon Klein's "Erlangen" approach (Klein 1872) through allowing for non-reversible transformations of spaces.

In fact Lawvere in his early papers (1964) and (1966) refers to the standard formal method and presents his proposed theories as first-order formal theories of categories. However in his later work (2003) speaking about category-theoretic foundations of mathematics Lawvere no longer refers to the formal method and opts for understanding of the notion of foundations "in a common-sense way rather than in the speculative way of the Bolzano-Frege-Peano-Russell tradition". I find this change of Lawvere's view significant and suggest to consider the categorical method of theory building pioneered by Lawvere as a genuine alternative to the formal method (rather than as something commonsensical). Remark that such categories as the category of (all) sets, all groups, all manifolds, etc. cannot possibly have isomorphic copies. For an isomorphic copy of, say, "the" category of sets cannot be anything else but another category of sets. But all the sets are supposed to be already there in the former category. As far as such a category is specified up to isomorphism it turns to be unique. From the formal point of view this is the most suspicious and makes one to suggest that notions of categories of "all" objects of a given type are illegitimate. This view is defended in (Eilenberg&MacLane1945) on a set-theoretic ground. However since then categories of this sort became ubiquitous in mathematical practice. In fact the talk of the category of "all" sets and other similar categories refers to categorical construals of corresponding concepts (of the concept of set and others) making the issue of "size" irrelevant. ("Local smallness" in categories is a different issue, which I leave aside.) So it is not surprising that the notion of isomorphism of categories turns to be, as Gelfand and Manin's put this, "useless" (2003, p.70) even if not illegitimate. This fact provides an additional support to the idea of developing categorical methods of theory-building independently of formal methods: unless a class of isomorphic models is available the very notion of formal looses its sense.

Before providing details about categorical methods of theory-building let's discuss what we may hope to get. A reasonable task is, for example, to identify a category as the category of Riemannian manifolds or as the category of Euclidean spaces. (It goes without saying that concepts of Riemannian manifold, Euclidean spaces and other mathematical concepts can be conceived in different ways. Speaking about "the" category of manifolds, etc. I mean the corresponding category of "all" objects of a given sort corresponding to the chosen version of a given concept. I don't mean to rule out the possibility that different categories can be called "categories of Riemannian manifold".)

This is what earlier in this paper I called a construal of a given concept "up to interpretation" or "up to arbitrary morphism". If we now compare this expected result with what one gets through application of formal methods we can notice two things. (i) Formal methods normally describe what in a corresponding categorical construction is a single object (up to isomorphism) like the Euclidean plane. In this sense categorical theories (now I mean theories built by category-theoretic methods) are more comprehensive than corresponding formal theories. (ii) At the same the category of Riemannian manifolds doesn't comprise all arithmetical and other "external" models of geometrical spaces allowed from the formal viewpoint. In a categorical framework such models can be construed as more specific constructions with functors to other categories. So the categorical approach viewed against the formal one not only brings new objects under the scope of the same conceptual scheme but also discriminates between models formally seen as essentially identical. Thus the categorical approach supports a more traditional view according to which, pairs of numbers may *represent* rather than *be* geometrical points. Remarkably this discrimination is made on a ground which a formalist would accept: the two kinds of things are distinguished by their different "behaviour" but not by anything like their "intrinsic nature"; although in some special contexts points and pairs of numbers may behave similarly in wider contexts taken into account by categorical methods these things behave quite differently.

Let me now point to a categorical method of theory-building suggested in (Lawvere 1963-2004) and called *functorial semantics*. It looks like a simulation of formal axiomatic method by categorical means but at the same time it demonstrates essential differences. Instead of writing axioms with usual strings of formulae one encodes axioms into a special "syntactic" category **T** that plays the role of "formal theory". Like in the case of standard semantics one assumes a "background" semantic category **B**, which is usually taken to be the category of sets but can be chosen differently. Models of **T** are functors of the form **f: T-->B**. One advantage of this construction is that it allows for different notions of model (only functors having certain properties count as models). Another remarkable fact is that under rather general conditions **T** can be embedded into a category **M(T,B)** of its functorial models. This definitely changes the whole idea of theory as a structure over and above all its possible models and suggests the view on a theory as "generic model" (Lawvere 1963-2004, p.19), which generates other models like circles and straight lines generate further constructions in Euclid's *Elements*. The functorial semantics also shows that the requirement of categoricity (in Veblen's sense) is as much unrealistic as unreasonable: although "good" categorical properties of **M(T,B)** are desirable there is no good reason to require that this category reduces to a single object.

Remark a different role of logic within formal and the categorical approaches. The formal approach fits the traditional view on logic (dating back to Aristotle) as the most general theory of reasoning, or more precisely, of the *correct* reasoning. In order to build a formal theory one "starts with" a system of formal logical calculus, and then writes down axioms introducing new non-logical symbols and thinking about their interpretations. So the idea of formal axiomatic method implies at least a weak version of logicism according to which any mathematical theory is based on (even if not reduces to) a certain system of logic. A categorical reconstruction of mathematical concepts unlike a formal reconstruction doesn't generally start with logic. It starts with the notion of category and proceeds with further categorical constructions. Permissible categorical constructions can be specified through constructive postulates. Consider this postulate, which I deliberately put in the form similar to Euclid's postulates:

*Given morphisms A-->B and B-->C to produce morphism A-->C*

A notion of logic can be then recaptured through associating a formal logical calculus with a category having appropriate "logical" properties (Bunge 1984, Makkai& Reyes 1977). This gives the notion of *internal logic* of a category. (As one might expect the internal logic of the category of sets is Classical logic.) So the categorical approach unlike the formal one is compatible with the view (held by otherwise so different thinkers as Boole, Poincaré and Brouwer) that logic assumes certain mathematical principles and so cannot be viewed as a basis for mathematics. This latter view on relationships between mathematics and logic better copes with the pluralism about logic (Beall&Restall 2000) than the logicist view mentioned above.

Lawvere's functorial semantics has been developed for a special case of algebraic theories and so it cannot be immediately used as a method of theory-building applicable in all areas of mathematics. Since then a lot of technical work has been done in related fields of categorical logic and categorical model theory. Nevertheless the categorical method of theory-building remains a work in progress and doesn't exist to the date in a standard form. For a historical introduction and further references I refer the reader to (Bell 2005).

3. Reversibility of mathematical reasoning and conceptual development of mathematics

We have seen that taking into consideration of non-reversible transformations between mathematical objects and treating these transformations on equal footing with isomorphisms has quite dramatic consequences for mathematics and its philosophy. Category theory is the general theory of non-reversible transformations (morphisms). Category theory meets the *hermeneutic challenge* of geometry of the end of 19-th century in a far more radical way than Hilbert's formal method: it makes morphisms, which can be viewed as *elementary interpretations*, into building blocks of mathematical constructions. (The notion of *object* of a category is in fact redundant: identity morphisms are sufficient to make the categorical machinery work properly.) In this section I shall try to outline some features of the new notion of mathematics brought about by this development.

Let me first show that considering of non-reversible morphisms on equal footing with isomorphisms is less innocent than it might seem. Consider this question: Is the operation of addition 7+5-->12 reversible or not? The question can be understood in different senses and so given different answers. The operation (+5) can be cancelled (reversed) by this subtraction: 12-5-->7, so it is reversible. However given that 12 is obtained as a sum of two summands there is no way to find out what these summands are. In this latter sense the operation (in fact a different operation) is non-reversible. But yet in a different (and perhaps not a "properly mathematical") sense the latter operation is still reversible: when 7 and 5 are summed up and bring 12 the summands 7, 5 don't perish but survive! This allows for writing the usual equality sign instead of arrow: 7+5=12. In *this* latter sense any mathematical operation and any categorical morphism is reversible. This suggests thinking of the notions of mathematical operation and transformation as mere metaphors describing particular *relations* between mathematical objects. In this latter view when 7 and 5 are summed up "nothing happens" indeed: the story of emergence of 12 out of 7 and 5 is just a way to say that the three numbers hold a particular ternary relation. This way to explain away references to processes and operations in mathematical discourses is known at least since Plato and in 20th century it has been made again popular by Frege, Russell and their followers. I shall not discuss this Platonic view systematically here but only remark that it ceases to be plausible as far as mathematics is viewed as a human activity going on in space and time. A computer performing the operation 7+5-->12 may keep or not keep the summands in its memory after the operation is done. The same is true for humans. One may argue after Plato that these facts have nothing do with numbers themselves but I think that the task of reconstruction of a notion of number from the relevant conceptual dynamics empirically observed in mathematical classrooms and elsewhere is interesting and anyway worth trying.

Observe that if during a mathematical reasoning one forgets where he or she has started from this certainly disqualifies the reasoning. So the reversibility is certainly required in this case. This basic reversibility of mathematical reasoning is made explicit in the *Elements* where each proposition is repeated twice: immediately before and immediately after its proof (but before *Q.E.D.*) Now granting this basic reversibility one might argue that mathematics is ultimately formal while non-reversible mathematical transformations are superfluous structures construed on a formal basis. However the argument can be met through taking a larger-scale conceptual dynamics into consideration. The larger-scale dynamics shows that the reversibility of mathematical reasoning is not so fundamental as it seems. At larger temporal (and perhaps also spatial) scales mathematics is obviously non-reversible. When Pythagorean theorem is taught in school a teacher may reasonably aim at reversibility of all interpretations of this theorem given by different pupils since this shows that all the pupils have learnt one and the same thing. (Interestingly, not *any* kind of reversibility is desired in this situation. The merely phonetic reversibility will not do: when all the pupils in a class utter the statement of the theorem and its proof in exactly the same words the teacher may suspect that none of them in fact *understands* it.) But this standard obviously doesn't apply at historical scales. One can hardly talk about a common *form* of the Pythagorean theorem invariantly preserved since Ancient times to today (except the form of the famous diagram known as *Bride's Chair*, *Pythagoras' Trousers* and a number of different names). This is moreover true for developed mathematical theories. The continuity of development of mathematics can be better described in terms of non-reversible interpretations between different fragments of mathematics, including historically remote fragments. Quasi-eternal concepts like *the* Pythagorean theorem (or, say, *natural number*) can be best understood as epiphenomena of a continuous non-reversible conceptual dynamics. As the example of the Pythagorean theorem clearly shows transformations involved into this dynamics don't reduce neither to isomorphisms (which is obvious for otherwise mathematics couldn't develop) nor to *monomorphisms* (embeddings) of older contents into new ones. Like any other science mathematics not only acquires new knowledge but also constantly revises its older contents and throws some of them away. The cumulative model of development of science and mathematics is oversimplified even if it allows for occasional "revolutions" (Kuhn 1962). In a categorical framework such oversimplified assumptions no longer look as "natural". A categorical analysis makes it clear that to keep a certain branch of science (mathematical or not) at the same fixed point of its development is not a trivial task (as anybody involved into the educational business knows by experience). I suggest that this task is in fact impossible: without new research science and mathematics rapidly corrupt but not just cease to develop.

Thus however fundamental might seem the reversibility of mathematical reasoning, and however plausible be the Platonic view on mathematical objects as eternal, a wider outlook reveals that in fact these phenomena are local. Assuming the notion of non-reversible morphism as basic a mathematician is no longer obliged to think of mathematics only in terms of its today's slice.

4. Conclusion: Mathematical Structuralism.

In the philosophy of mathematics of the 20-th century *Structuralism* has been first associated with Bourbaki's fundamental *Les Eléments des mathématiques* aiming at reconstruction of mathematics in set-theoretic terms. *Les Eléments* starts with a version of Set theory, and then introduce various mathematical concepts as "structured sets". Take a set *G* and associate with any ordered pair of its elements a third element of the same set;

an obvious axiomatic description of this construction turns it into a *group*. In the given context the set $G$ can be called a "set equipped with the group structure". The idea is that the "structure" is put here on the top of an arbitrarily chosen "background set". Topological spaces, rings, modules and many other (arguably all) mathematical concepts can be construed similarly. A large part of mathematics of the second part of the 20-th century has been made in this framework. Bourbaki's colloquial metaphysics resembles Aristotle's: mathematical objects are "made of" the same set-theoretic "matter" and distinguished by their "forms" (structures).

Chapter 4 of (Bourbaki 1968) presents a general theory of structures. Here Bourbaki distinguishes basic kinds of structures like the *final* and the *initial* structures, which don't correspond to popular mathematical concepts like that of group or topological space. A reader familiar with the Category theory can easily identify them with corresponding category-theoretic concepts. So Category theory suggests to improve upon Bourbaki's colloquial metaphysics and think of "forms" or "structures" in the abstraction of the set-theoretic "matter"! Hence the claim that the Category theory provides a further support for the Mathematical Structuralism (Awodey 1996, 2004; MacLane 1996), which is the view that only structures (but not the "matter") count.

Here is a recent official definition of mathematical structuralism (Hellman, forthcoming):

"Structuralism is a view about the subject matter of mathematics according to which what matters are structural relationships in abstraction from the intrinsic nature of related objects."

Like in the case of Hilbert's *Grundlagen* I suggest to distinguish carefully between the two following features implied by this definition. The first is epistemological (and possibly also) ontological *holism* about mathematical matters: no mathematical object can be thought of outside of its relations with other mathematical objects. The second feature, which I shall call exchangism, is more specific. It is the assumption of possibility of exchanging the relata of a given relation by some other relata, so that the "structural relationships" or the "structure" remains invariant. Consider group $G$ construed as just described and replace each element of its underlying set by another element of another set. This operation obviously doesn't change much: what we get is just another isomorphic copy $G'$ of the "same" group. Similarly any model of a given formal theory can be replaced by another isomorphic model. Structuralism says that all mathematical theories and concepts function in this way. Hellman quite rightly traces this view back to Hilbert's *Grundlagen*.

Let's now test the two features of Structuralism against Category theory and the prospective categorical method of theory-building suggested above. Categorical approach certainly pushes mathematical holism further forward. While the formal approach suggests to chose a favourite model of a given formal theory (but keep in mind that this model can be exchanged for another isomorphic model) the categorical approach suggests working with a whole category of models at once. But the "exchangism" generally doesn't apply to Category theory and to prospective categorical mathematics. To see this remark that "exchanges" are reversible transformations (think about $G$ and $G'$) while categorical morphisms, generally speaking, are not. As I have already argued above, the idea that given, say, the category of sets and functions one may exchange sets and functions for groups and group homomorphisms, then again for topological spaces and continuous transformation, etc., and finally to conceive of an abstract notion of category as a structure invariant through such exchanges is plainly wrong. For categories just mentioned are not isomorphic (functors between them are not reversible), so there is nothing invariant here. So Category theory and categorical mathematics doesn't support Mathematical Structuralism or at least the version of Structuralism explicated above. Remarkably the categorical approach defended here doesn't rely on anything like "intrinsic nature" of mathematical objects either. Moreover it apparently squares well with this basic idea behind the structuralist approach: mathematics studies relations between things but not things themselves. The notion of mathematical structure has been developed as an attempt of thinking "relations without relata" or more precisely "relations indifferent to their relata". Considering categorical morphisms as relations (in a broad philosophical rather than technical sense) one may argue that Category theory provides even better mathematical realisation of this project than Bourbaki's set-theoretic machinery. Given that objects in a category are identified with their identity morphisms the category is construed out of morphisms ("relations") and nothing else.

The idea of categorical mathematics, which I defend in this paper, is indeed a further development of the idea of structuralist mathematics. But in order to make the difference clear I choose not to qualify the categorical approach as structuralist. The talk of "structural relationships in abstraction from the intrinsic nature of related objects" clearly implies the notion of structure as a form-concept. Trying to generalise this notion to the effect of making it into a category-concept one looses anything like its usual understanding. Given a class $C$ of isomorphic mathematical objects one may always associate with it a new abstract object $S$ conceived as "structure" shared by all members of $C$; then $S$ can be thought of " in abstraction from the intrinsic nature" of members of $C$. This Fregean notion of abstraction seems to be relevant to Hellman's definition of Structuralism. But I don't see how this notion of abstraction can be reasonably generalised in order to allow for the case of non-reversible morphisms. Given an isomorphism $A<-->B$ one may identify $A, B$ "up to isomorphism" and conceive the result of this identification as a new abstract object $C$. However given a non-reversible morphism $A-->B$ one cannot perform anything similar. So the notion of abstract structure stipulated over and above all its isomorphic instantiations doesn't survive the suggested generalisation.

The idea of categorical mathematics pushes further forward the anti-essentialist stance of Structuralism. From a categorical viewpoint the notion of mathematical structure looks like just another kind of "mathematical essence". I can see that structures played and continue to play an important role in mathematics but I claim that the structuralist view on mathematics as a general science of structures is not adequate to today's mathematics and can hardly motivate new mathematical developments. The view on mathematics as a general science of categories is much broader (since the notion of category is far more general than the notion of structure) and provides mathematicians and philosophers with a new perspective yet to be explored.